\input amstex
\input amsppt.sty
\magnification=\magstep1
\hsize=33.5truecc
\vsize=23truecm
\baselineskip=16truept
\NoBlackBoxes
\TagsOnRight \pageno=1 \nologo
\def\Z{\Bbb Z}
\def\N{\Bbb N}

\def\l{\left}
\def\r{\right}
\def\bg{\bigg}
\def\({\bg(}
\def\[{\bg\lfloor}
\def\){\bg)}
\def\]{\bg\rfloor}
\def\t{\text}
\def\f{\frac}

\def\bi{\binom}
\def\eq{\equiv}

\def\ls{\leqslant}
\def\gs{\geqslant}
\def\mo{\roman{mod}}

\def\da{\delta}

\def\M#1#2{\thickfracwithdelims[]\thickness0{#1}{#2}_q}

\def\Proof{\noindent{\it Proof}}

\def\Remark{\medskip\noindent{\it  Remark}}

\def\Ack{\medskip\noindent {\bf Acknowledgments}}
\hbox {Adv. Appl. Math. 136 (2022), Article ID 102319.}
\bigskip
\topmatter
\title On Motzkin numbers and central trinomial coefficients\endtitle
\author Zhi-Wei Sun\endauthor
\leftheadtext{Zhi-Wei Sun}
\rightheadtext{On Motzkin numbers and central trinomial coefficients}
\affil Department of Mathematics, Nanjing University\\
 Nanjing 210093, People's Republic of China
  \\  zwsun\@nju.edu.cn
  \\ {\tt http://maths.nju.edu.cn/$\sim$zwsun}
\endaffil
\abstract The Motzkin numbers $M_n=\sum_{k=0}^n\bi n{2k}\bi{2k}k/(k+1)$ $(n=0,1,2,\ldots)$
and the central trinomial coefficients $T_n$ ($n=0,1,2,\ldots)$ given by the constant term of $(1+x+x^{-1})^n$, have many combinatorial interpretations. In this paper we establish
the following surprising arithmetic properties
of them with $n$ any positive integer:
$$\f2n\sum_{k=1}^n(2k+1)M_k^2\in\Z,$$
$$\f{n^2(n^2-1)}6\,\bigg|\,\sum_{k=0}^{n-1}k(k+1)(8k+9)T_kT_{k+1},$$
and also
$$\sum_{k=0}^{n-1}(k+1)(k+2)(2k+3)M_k^23^{n-1-k}=n(n+1)(n+2)M_nM_{n-1}.$$
\endabstract
\thanks 2020 {\it Mathematics Subject Classification}. \,Primary 05A10, 05A19;
Secondary  11A07, 11B75.
\newline\indent Supported by the National Natural Science
Foundation of China (grant 11971222).
\endthanks

\endtopmatter
\document

\heading{1. Introduction}\endheading

In combinatorics, the Motzkin number $M_n$ with $n\in\N=\{0,1,2,\ldots\}$
is the number of lattice paths from the point $(0,0)$ to the point
$(n,0)$ which never dip below the line $y=0$ and are made up only of the allowed steps $(1,0)$ (east), $(1,1)$ (northeast) and $(1,-1)$ (southeast).
It is well known that
$$M_n=\sum_{k=0}^{\lfloor n/2\rfloor}\bi n{2k}C_k$$
where $C_k$ denotes the Catalan number $\bi{2k}k-\bi{2k}{k+1}=\bi{2k}k/(k+1)$.

For $n\in\N$, the central trinomial coefficient $T_n$ is the constant term in the expansion of $(1+x+x^{-1})^n$.
By the multi-nomial theorem, we see that
$$T_n=\sum_{k=0}^{\lfloor n/2\rfloor}\bi n{2k}\bi{2k}k=\sum_{k=0}^n\bi nk\bi{n-k}k.$$
It is known that $T_n$ coincides with the number of lattice paths from the point $(0, 0)$ to $(n, 0)$
with only allowed steps $(1,0)$ (east), $(1, 1)$ (northeast) and $(1, -1)$ (southeast).

The Motzkin numbers, the Catalan numbers and the central trinomial coefficients arise naturally in enumerative combinatorics. As the Fibonacci numbers arising from combinatorics have rich number-theoretic properties, we think that important combinatorial quantities like
$M_n$ and $T_n$ with $n\in\N$ should also have nice arithmetic properties. For example, in [S14a] we conjectured that for any $n\in\Z^+=\{1,2,3,\ldots\}$
the arithmetic mean of the $n$ numbers $(8k+5)T_k^2\ (k=0,\ldots,n-1)$ is always an integer,
and this was later confirmed by Y.-P. Mu and the author [MS] via symbolic computation.
Motivated by congruence properties of such numbers, we found in [S14b, S20] many series for $1/\pi$
involving central trinomial coefficients or their extensions.
For example, in [S20, Section 10] we conjectured the combinatorial identity
$$\sum_{k=1}^\infty\f{(105k-44)T_{k-1}}{k^2\bi{2k}k^23^{k-1}}=\f{5\pi}{\sqrt3}+6\log3$$
based on the conjectural congruence
$$p^2\sum_{k=1}^{p-1}\f{(105k-44)T_{k-1}}{k^2\bi{2k}k^23^{k-1}}
\eq11\l(\f p3\r)+\f p2\l(13-35\l(\f p3\r)\r)\pmod{p^2},$$
where $p$ is a prime greater than $3$ and $(-)$ is the Legendre symbol.
Thus it is interesting to investigate congruence properties of combinatorial quantities
like $M_n$ and $T_n$ with $n\in\N$, and the study in turn may stimulate us to find some new combinatorial identities.

Let $p>3$ be a prime. In [S14a, Conjecture 1.1(ii)] we conjectured
$$\sum_{k=0}^{p-1}M_k^2\eq(2-6p)\l(\f p3\r)\pmod{p^2},\ \ \sum_{k=0}^{p-1}kM_k^2\eq(9p-1)\l(\f p3\r)\pmod{p^2},$$
and
$$\sum_{k=0}^{p-1}T_kM_k\eq\f 43\l(\f p3\r)+\f p6\l(1-9\l(\f p3\r)\r)\pmod{p^2}.$$
The three supercongruences look curious and challenging.

Motivated by the above conjectures, we establish the following new results.

\proclaim{Theorem 1.1} {\rm (i) For any $n\in\Z^+$, we have
$$s(n):=\f2n\sum_{k=1}^n(2k+1)M_k^2\in\Z.\tag1.1$$

{\rm (ii)} For any prime $p>3$, we have
$$\sum_{k=0}^{p-1}(2k+1)M_k^2\eq12p\l(\f p3\r)\pmod{p^2}.\tag1.2$$
\endproclaim
\Remark\ 1.1. The values of $s(1),\ldots,s(10)$ are as follows:
$$6,\ 23,\ 90,\ 432,\ 2286,\ 13176,\ 80418,\ 513764,\ 3400518,\ 23167311.$$

\proclaim{Theorem 1.2} For any integer $n\gs2$, we have
$$\f{n^2(n^2-1)}6\,\bigg|\,\sum_{k=0}^{n-1}k(k+1)(8k+9)T_kT_{k+1}.\tag1.3$$
\endproclaim
\Remark\ 1.2. If we define
$$t(n):=\f 6{n^2(n^2-1)}\sum_{k=0}^{n-1}k(k+1)(8k+9)T_kT_{k+1}\ \ \ (n=2,3,\ldots),$$
then the values of $t(2),t(3),\ldots,t(10)$ are as follows:
$$51,\ 271,\ 1398,\ 8505,\ 54387,\ 367551,\ 2570931,\ 18510739,\ 136282347.$$

Let $b,c\in\Z$ and $n\in\N$. The generalized central trinomial coefficient $T_n(b,c)$
denotes the coefficient of $x^n$ in the expansion of $(x^2+bx+c)^n$ (cf. [S14a] and [S14b]). By the multi-nomial theorem, we see that
$$T_n(b,c)=\sum_{k=0}^{\lfloor n/2\rfloor}\bi n{2k}\bi{2k}kb^{n-2k}c^k.$$
The generalized Motzkin number $M_n(b,c)$ introduced in [S14a] is given by
$$M_n(b,c)=\sum_{k=0}^{\lfloor n/2\rfloor}\bi n{2k}C_kb^{n-2k}c^k.$$
Note that $T_n(1,1)=T_n$, $M_n(1,1)=M_n$, $T_n(2,1)=\bi{2n}n$ and $M_n(2,1)=C_{n+1}$.
Also, $T_n(3,2)$ coincides with the (central) Delannoy number
$$D_n=\sum_{k=0}^n\bi nk\bi{n+k}k=\sum_{k=0}^n\bi{n+k}{2k}\bi{2k}k,$$
 which counts lattice paths from $(0,0)$ to $(n,n)$ in which only east $(1, 0)$, north $(0, 1)$, and northeast $(1, 1)$ steps are allowed
(cf. R. P. Stanley [St99, p.\,185]). And $M_n(3,2)$ equals the little Schr\"oder number
$$s_{n+1}=\sum_{k=1}^{n+1}N(n+1,k)2^{n+1-k}$$
with the Narayana number $N(m,k)$ ($m\gs k\gs1$) given by
$$N(m,k):=\f1m\bi mk\bi m{k-1}\in\Z.$$
The little Schr\"oder numbers and the Narayana numbers also
have many combinatorial interpretations (cf. [St97] and [Gr, pp.\,268--281]).
See also [S11, S18b] for some congruences involving the Delannoy numbers or the little Schr\"oder numbers.

\proclaim{Theorem 1.3} Let $b,c\in\Z$ with $b\not=0$ and $d=b^2-4c\not=0$, and let $n\in\Z^+$. Then
$$b\f{n(n+1)}2\,\bigg|\,\sum_{k=1}^n kT_k(b,c)T_{k-1}(b,c)d^{n-k}\tag1.4$$
and
$$b\f{n^2(n+1)^2}4\,\bigg|\,3\sum_{k=1}^nk^3T_k(b,c)T_{k-1}(b,c)d^{n-k}.\tag1.5$$
Also,
$$\f{(2,n)}{n(n+1)(n+2)}\sum_{k=0}^{n-1}(k+1)(k+2)(2k+3)M_k(b,c)^2d^{n-1-k}\in\Z\tag1.6$$
and
$$\sum_{k=0}^{n-1}\f{(k+1)(k+2)(2k+3)}{n(n+1)(n+2)}M_k(b,c)^2(-d)^{n-1-k}
=\f{M_n(b,c)M_{n-1}(b,c)}b\in\Z,\tag1.7$$
where $(m,n)$ denotes the greatest common divisor of two integers $m$ and $n$.
\endproclaim
\Remark\ 1.3. For each $n\in\Z^+$, (1.7) with $b=c=1$ gives the curious identity
$$\sum_{k=0}^{n-1}(k+1)(k+2)(2k+3)M_k^23^{n-1-k}=n(n+1)(n+2)M_nM_{n-1}.\tag1.8$$

In the case $b=3$ and $c=2$, Theorem 1.3 yields the following consequence.

\proclaim{Corollary 1.1} For any $n\in\Z^+$ we have
$$3\f{n(n+1)}2\,\bigg|\,\sum_{k=1}^n kD_kD_{k-1},\ \f{n^2(n+1)^2}4\,\bigg|\,\sum_{k=1}^n k^3D_kD_{k-1},\tag1.9$$
$$\f{n(n+1)(n+2)}{(2,n)}\,\bigg|\,\sum_{k=1}^{n}k(k+1)(2k+1)s_k^2,\tag1.10$$
and
$$\f1{n(n+1)(n+2)}\sum_{k=1}^nk(k+1)(2k+1)(-1)^{n-k}s_k^2=\f{s_ns_{n+1}}3\in\Z.\tag1.11$$
\endproclaim

Theorems 1.1-1.3 are quite sophisticated and their proofs need various techniques.
We will prove Theorems 1.1-1.3 in Sections 2-4 respectively. In Section 5 we are going to pose some related conjectures for further research.

\heading{2. Proof of Theorem 1.1}\endheading

For $n\in\Z^+$, in [S18b] we introduced the polynomial
$$s_n(x):=\sum_{k=1}^nN(n,k)x^{k-1}(x+1)^{n-k}\tag2.1$$
for which $s_n(1)$ is just the little Schr\'oder number $s_n$.
For $n\in\N$, define
$$S_n(x)=\sum_{k=0}^n\bi nk\bi{n+k}k\f{x^k}{k+1}=\sum_{k=0}^n\bi {n+k}{2k}C_kx^k.\tag2.2$$
Then $S_n(1)$ equals the large Schr\"oder number $S_n$
which counts the lattice paths from the point $(0,0)$ to $(n,n)$ with steps
$(1,0),(0,1)$ and $(1,1)$ that never rise above the line $y=x$.
As proved in [S18b], we have
$$S_n(x)=(x+1)s_n(x)\ \quad\t{for all}\ n\in\Z^+.\tag2.3$$

\proclaim{Lemma 2.1}
{\rm (i)} For any $n\in\Z^+$ we have
$$n(n+1)s_n(x)^2=\sum_{k=1}^n\bi{n+k}{2k}\bi{2k}k\bi{2k}{k+1}(x(x+1))^{k-1}.\tag2.4$$

{\rm (ii)} Let $b,c\in\Z$ with $d=b^2-4c\not=0$. For any $n\in\N$ we have
$$M_n(b,c)=(\sqrt d)^ns_{n+1}\l(\f{b/\sqrt d-1}2\r).\tag2.5$$
\endproclaim
\Proof. As $(x+1)s_n(x)=S_n(x)$ by (2.3), the identity (2.4) has the equivalent version
$$n(n+1)S_n(x)^2=\sum_{k=1}^n\bi{n+k}{2k}\bi{2k}k\bi{2k}{k+1}x^{k-1}(x+1)^{k+1}$$
which appeared as [S12a, (2.1)]. So (2.4) holds. The identity (2.5) was proved in [S18b, Lemma 3.1]. \qed

\Remark\ 2.1. For $n\in\N$ and $b,c\in\Z$ with $b^2\not=4c$, by combining the two parts of Lemma 2.1 we obtain that
$$M_n(b,c)^2=\f1{(n+1)(n+2)}\sum_{k=1}^{n+1}\bi{n+k+1}{2k}\bi{2k}k\bi{2k}{k+1}c^{k-1}(b^2-4c)^{n+1-k}.\tag2.6$$

\proclaim{Lemma 2.2} For any $n\in\Z^+$ we have
$$\aligned&\sum_{k=1}^n(2k+1)M_k^2
\\=&\sum_{k=0}^{n+1}\f{(4n-2k+3)(n+k+2)}{n+2}\bi{n+k+1}{2k}\bi{2k}{k}\bi{2k+1}{k}(-3)^{n+1-k}.
\endaligned\tag2.7$$
\endproclaim
\Proof. In view of (2.6), we have
$$\align \sum_{k=0}^n(2k+1)M_k^2=&\sum_{k=0}^n\f{2k+1}{(k+1)(k+2)}\sum_{j=1}^{k+1}\bi{k+j+1}{2j}\bi{2j}j\bi{2j}{j+1}(-3)^{k+1-j}
\\=&\sum_{k=0}^n\f{2k+1}{(k+1)(k+2)}\sum_{l=0}^{k}\bi{k+l+2}{2l+2}\bi{2l+2}{l+1}\bi{2l+2}{l}(-3)^{k-l}
\\=&\sum_{k=0}^n\sum_{l=0}^nF(k,l),
\endalign$$
where
$$F(k,l):=\f{2k+1}{(k+1)(k+2)}\bi{k+l+2}{2l+2}\bi{2l+2}{l+1}\bi{2l+2}{l}(-3)^{k-l}.$$
By the telescoping method developed by Chen, Hou and Mu [CHM] and applied by
Mu and Sun [MS], the double sum can be reduced to a single sum:
$$\sum_{k=0}^n\sum_{l=0}^nF(k,l)=1+(4n+3)(-3)^{n+1}+\sum_{j=0}^n(-3)^{n-j}\f{(4n-2j+1)(n+j+3)!(2j+3)!}{(n+2)(n-j)!(j+2)(j+1)!^4}.\tag2.8$$
Therefore
$$\align&\sum_{k=1}^n(2k+1)M_k^2
\\=&\sum_{j=-1}^n(-3)^{n-j}\f{(4n-2j+1)(n+j+3)!(2j+3)!}{(n+2)(n-j)!(j+2)(j+1)!^4}
\\=&\sum_{k=0}^{n+1}(-3)^{n+1-k}\f{(4n-2k+3)(n+k+2)!(2k+1)!}{(n+2)(n+1-k)!(k+1)k!^4}
\\=&\sum_{k=0}^{n+1}\f{(4n-2k+3)(n+k+2)}{n+2}\bi{n+k+1}{2k}\bi{2k}{k}\bi{2k+1}{k}(-3)^{n+1-k}
\endalign$$
and this concludes the proof. \qed

For each integer $n$ we set
$$[n]_q=\f{1-q^n}{1-q},$$
which is the usual $q$-analogue of $n$. For any $n\in\Z$, we define
$$\M n0=1\quad\t{and}\quad \M nk=\f{\prod_{j=0}^{k-1}[n-j]_q}{\prod_{j=1}^k[j]_q}\ \ \t{for}\ k=1,2,3,\ldots.$$
Obviously $\lim_{q\to1}\M nk=\bi nk$ for all $k\in\N$ and $n\in\Z$. It is easy to see that
$$\M nk=q^k\M{n-1}k+\M{n-1}{k-1}\quad \ \t{for all}\ k,n=1,2,3,\ldots.$$
By this recursion, $\M nk\in\Z[q]$ for all $k,n\in\N$. For any integers $a,\,b$ and $n>0$, clearly
$$a\eq b\pmod n\ \Longrightarrow\ [a]_q\eq[b]_q\pmod{[n]_q}.$$

Let $n$ be a positive integer. The cyclotomic polynomial
$$\Phi_n(q):=\prod^n\Sb a=1\\(a,n)=1\endSb\l(q-e^{2\pi ia/n}\r)\in\Z[q]$$
is irreducible in the ring $\Z[q]$. It is well-known that
$$q^n-1=\prod_{d\mid n}\Phi_d(q).$$
Note that $\Phi_1(q)=q-1$.

\proclaim{Lemma 2.3} For any $a,b\in\N$ and $n\in\Z^+$, we have
$$\sum_{k=0}^{n-1}\M{n+1}k^a\M{n+k}k^b\M{2k}k[k+2]_q(-[3]_q)^{n-1-k}\eq0\pmod{[n]_q}.\tag2.9$$
\endproclaim
\Proof. (2.9) is trivial in the case $n=1$. Below we assume $n>1$. As
$$[n]_q=\prod_{1<d\mid n}\Phi_d(q)$$
and $\Phi_2(q),\Phi_3(q),\ldots$ are pairwise coprime, it suffices to show that the sum in (2.9)
is divisible by $\Phi_d(q)$ for any given divisor $d>1$ of $n$.

 A well-known $q$-Lucas theorem (see, e.g., [O]) states that if $a,b,d,s,t\in\N$ with $s<d$ and $t<d$ then
$$\M{ad+s}{bd+t}\eq\bi ab\M st\pmod{\Phi_d(q)}.$$
Let $S$ denote the sum in (2.9) and write $n=dm$ with $m\in\Z^+$. Then
$$\align S=&\sum_{j=0}^{m-1}\sum_{r=0}^{d-1}\M{md+1}{jd+r}^a\M{md+jd+r}{jd+r}^b\M{2jd+2r}{jd+r}[jd+r+2]_q(-[3]_q)^{md-1-(jd+r)}
\\\eq&\sum_{j=0}^{m-1}\sum_{r=0}^{d-1}\bi mj^a\M 1r^a\bi{m+j}j^b\M{r}{r}^b\M{2jd+2r}{jd+r}[r+2]_q(-[3]_q)^{(m-j)d-(r+1)}
\\\eq&\sum_{j=0}^{m-1}\bi mj^a\bi{m+j}j^b\sum_{r=0}^1\M{2jd+2r}{jd+r}[r+2]_q(-[3]_q)^{(m-j)d-(r+1)}
\\\eq&\sum_{j=0}^{m-1}\bi mj^a\bi{m+j}j^b\bi{2j}j\M{0}{0}[2]_q(-[3]_q)^{(m-j)d-1}
\\&+\sum_{j=0}^{m-1}\bi mj^a\bi{m+j}j^b[1+2]_q(-[3]_q)^{(m-j)d-2}\times\cases\bi{2j+1}j\M{0}{1}&\t{if}\ d=2,
\\\bi{2j}j\M{2}1&\t{if}\ d>2,\endcases
\\\eq&0\pmod{\Phi_d(q)}.
\endalign$$
(Note that $[2]_q=1+q=\Phi_2(q)$.)
This concludes the proof. \qed

\proclaim{Lemma 2.4} For any prime $p>3$ we have
$$\sum_{k=1}^{p-1}\f{\bi{2k}k}{k3^k}\eq\f{3^{p-1}-1}p\pmod p.\tag2.10$$
\endproclaim
\Proof. Let $u_n=(\f n3)$ for $n\in\N$. Then
$u_0=0$, $u_1=1$ and $u_{n+1}=-u_n-u_{n-1}$ for all $n=1,2,3,\ldots$. Applying [S12b, Lemma 3.5] with $m=1$, we obtain
$$\sum_{k=1}^{p-1}\f{\bi{2k}k}{k3^k}\eq\f{(-3)^{p-1}-1}p-\f12\l(\f{-3}p\r)\f{u_{p-(\f{-3}p)}}p\pmod p.$$
Note that $u_{p-(\f{-3}p)}=0$ since $p\eq (\f{-3}p)\pmod3$. So (2.10) holds. \qed

\medskip
\noindent{\it Proof of Theorem} 1.1. (i) Observe that
$$\f 4{n+2}\eq\cases 4/2=2\pmod n&\t{if}\ 2\nmid n,
\\2/(n/2+1)\eq2\pmod n&\t{if}\ 2\mid n.\endcases$$
Thus, for each $k\in\{1,\ldots,n+1\}$, we have
$$2\times\f{\bi{2k}k}{n+2}=\f 4{n+2}\bi{2k-1}k\eq2\bi{2k-1}k=\bi{2k}k\pmod n.$$
Combining this with (2.7) we see that
$$\align &2\sum_{k=1}^n(2k+1)M_k^2
\\\eq&2(4n+3)(-3)^{n+1}
\\&+\sum_{k=1}^{n+1}(4n-2k+3)(n+k+2)\bi{n+k+1}{2k}\bi{2k}k\bi{2k+1}k(-3)^{n+1-k}
\\\eq&-\sum_{k=0}^{n+1}(2k-3)(k+2)\bi{n+k+1}{n+1}\bi{n+1}k\bi{2k+1}k(-3)^{n+1-k}
\\\eq&-\sum_{k=0}^{n+1}(2k-3)(k+2)\f{n+k+1}{n+1}\bi{n+k}k\bi{n+1}k(2k+1)C_k(-3)^{n+1-k}
\\\eq&-\sum_{k=0}^{n+1}(2k-3)(k+2)(k+1)\bi{n+k}k\bi{n+1}k(2k+1)C_k(-3)^{n+1-k}\pmod{n}.
\endalign$$
For each $k=0,\ldots,n+1$, clearly
$$k(k-1)\bi{n+1}k=n(n+1)\bi{n-1}{n+1-k}\eq0\pmod n.$$
Since $(2k-3)(2k+1)=4k(k-1)-3$, by the above we have
$$2\sum_{k=1}^n(2k+1)M_k^2
\eq-\sum_{k=0}^{n+1}\bi{n+1}k\bi{n+k}k\bi{2k}k(k+2)(-3)^{n+2-k}\pmod n.$$
Note that
$$\align&\sum_{k=n}^{n+1}\bi{n+1}k\bi{n+k}k\bi{2k}k(k+2)(-3)^{n+2-k}
\\=&\bi{n+1}n\bi{2n}n^2(n+2)(-3)^2+\bi{2n+1}{n+1}\bi{2n+2}{n+1}(n+3)(-3)
\\\eq&18\bi{2n}n^2-18\(\f{2n+1}{n+1}\bi{2n}n\)^2\eq0\pmod{n}.
\endalign$$
Therefore
$$2\sum_{k=1}^n(2k+1)M_k^2
\eq27\sum_{k=0}^{n-1}\bi{n+1}k\bi{n+k}k\bi{2k}k(k+2)(-3)^{n-1-k}\ (\mo\ n).\tag2.11$$

By (2.9) with $a=b=1$ and $q=1$, we have
$$\sum_{k=0}^{n-1}\bi{n+1}k\bi{n+k}k\bi{2k}k(k+2)(-3)^{n-1-k}\eq0\pmod n.$$
Combining this with (2.11) we immediately obtain the desired (1.1).

(ii) Applying (2.7) with $n=p-1$, we get
$$\align \sum_{k=1}^{p-1}(2k+1)M_k^2
=&\sum_{k=0}^p\f{(4p-2k-1)(p+k+1)}{p+1}\bi{p+k}{2k}\bi{2k}k\bi{2k+1}k(-3)^{p-k}
\\=&\sum_{k=1}^{p-1}\f{(4p-2k-1)(p+k+1)}{p+1}\bi{p}k\bi{p+k}k\f{2k+1}{k+1}\bi{2k}k(-3)^{p-k}
\\&+(4p-1)(-3)^p+\f{(2p-1)(2p+1)}{p+1}\bi{2p}p\f{2p+1}{p+1}\bi{2p}p
\\\eq&3\sum_{k=1}^{p-1}\f pk\bi{p-1}{k-1}(2k+1)^2\f{\bi{2k}k}{(-3)^k}+(3-12p)3^{p-1}-\l(2\bi{2p-1}{p-1}\r)^2
\\\eq&-3p\sum_{k=1}^{p-1}\l(4k+4+\f1k\r)\f{\bi{2k}k}{3^k}+3^p-12p-4\pmod{p^2}
\endalign$$
with the aid of Wolstenholme's congruence $\bi{2p-1}{p-1}\eq1\pmod{p^3}$ (cf. [W]).
Combining this with (2.10) and noting that $\bi{-1/2}k=\bi{2k}k/(-4)^k$ for $k\in\N$, we obtain
$$\align-\f1{12p}\sum_{k=0}^{p-1}(2k+1)M_k^2
\eq&1+\sum_{k=1}^{p-1}\l(k+1\r)\bi{-1/2}k\f{(-4)^k}{3^k}
\\\eq&\sum_{k=0}^{(p-1)/2}\bi{(p-1)/2}k\l(-\f 43\r)^k+\sum_{k=1}^{(p-1)/2}k\bi{(p-1)/2}k\l(-\f43\r)^k
\\\eq&\l(1-\f43\r)^{(p-1)/2}-\f43\cdot\f{p-1}2\sum_{k=1}^{(p-1)/2}\bi{(p-3)/2}{k-1}\l(-\f43\r)^{k-1}
\\\eq&\l(\f{-3}p\r)+\f23\l(1-\f43\r)^{(p-3)/2}
\\\eq&\l(\f{-3}p\r)-2\l(\f{-3}p\r)=-\l(\f p3\r)\pmod p.
\endalign$$
This proves (1.2).

The proof of Theorem 1.1 is now complete. \qed

\heading{3. Proof of Theorem 1.2}\endheading

\proclaim{Lemma 3.1} Let $b,c\in\Z$ and $d=b^2-4c$. Then
$$b\sum_{k=0}^{n-1}(2k+1)T_k(b,c)^2(-d)^{n-1-k}=nT_n(b,c)T_{n-1}(b,c)\ \ \ \t{for any}\ n\in\Z^+,\tag3.1$$
and $$T_k(b,c)^2=\sum_{j=0}^k\bi{k+j}{2j}\bi{2j}j^2c^jd^{k-j}\quad\t{for all}\ k\in\N.\tag3.2$$
\endproclaim
\Remark\ 3.1. For (3.1) and (3.2), see [S14a, (1.19) and (4.1)].

\proclaim{Lemma 3.2} For any $n\in\Z^+$, we have
$$\sum_{k=0}^{n-1}k(k+1)(8k+9)T_kT_{k+1}=\f{(-1)^nn}6\sum_{k=0}^{n-1}\bi {n-1}k\bi{-n-1}kC_k3^{n-1-k}a(n,k),\tag3.3$$
where
$$a(n,k)=4k^2n^2-8kn^3-14k^2n-14kn^2-4n^3+13k^2-11kn-26n^2+39k+4n+26.$$
\endproclaim
\Proof. In light of (3.1) with $b=c=1$,
$$\align &\sum_{k=0}^{n-1}k(k+1)(8k+9)T_kT_{k+1}
\\=&\sum_{k=1}^nk(k-1)(8k+1)T_kT_{k-1}
\\=&\sum_{k=1}^n(k-1)(8k+1)\sum_{j=0}^{k-1}(2j+1)T_j^23^{k-1-j}
\\=&\sum_{j=0}^{n-1}(2j+1)T_j^2\sum_{k=j+1}^n(k-1)(8k+1)3^{k-1-j}.
\endalign$$
By induction, for each $j\in\N$ we have
$$\sum_{k=j+1}^m(k-1)(8k+1)3^{k-1-j}=\f14\l(3^{m-j}(16m^2-30m+21)-(16j^2-30j+21)\r)$$
for all $m=j+1,j+2,\ldots$. Thus, in view of the above and (3.2) with $b=c=1$, we get
$$\align &4\sum_{k=0}^{n-1}k(k+1)(8k+9)T_kT_{k+1}
\\=&\sum_{k=0}^{n}(2k+1)T_k^2\l(3^{n-k}(16n^2-30n+21)-(16k^2-30k+21)\r)
=\sum_{k=0}^n\sum_{l=0}^nF(k,l),
\endalign$$
where $F(k,l)$ denotes
$$(2k+1)\bi{k+l}{2l}\bi{2l}l^2(-3)^{k-l}\l(3^{n-k}(16n^2-30n+21)-(16k^2-30k+21)\r).$$
Via the telescoping method stated in [CHM, MS], the double sum can be reduced to a single sum:
$$\sum_{k=0}^n\sum_{l=0}^nF(k,l)=\f29\sum_{k=0}^{n-1}\f{a(n,k)(-3)^{n-k}(n+k)!(2k)!}{(n-k-1)!k!^4(k+1)}.\tag3.4$$
Therefore
$$\align&\sum_{k=0}^{n-1}k(k+1)(8k+9)T_kT_{k+1}
\\=&\f1{18}\sum_{k=0}^{n-1}\bi n{k+1}\bi{n+k}k\bi{2k}k(-3)^{n-k}a(n,k)
\\=&\f{(-1)^n}6\sum_{k=0}^{n-1}\f n{k+1}\bi {n-1}{k}\bi{-n-1}k\bi{2k}k3^{n-1-k}a(n,k)
\endalign$$
and hence (3.3) holds. \qed

\proclaim{Lemma 3.3} For any $n\in\Z^+$, we have
$$n^2-1\ \bigg|\ \sum_{k=0}^{n-1}\bi{n-1}{k}\bi{-n-1}kC_k3^{n-1-k}a(n,k)\tag3.5$$
with $a(n,k)$ given in Lemma 3.2.
\endproclaim
\Proof. It suffices to show that $n^2-1$ divides $\bi{n-1}k\bi{-n-1}ka(n,k)$ for any fixed $k\in\{0,\ldots,n-1\}$.
Clearly,
$$\align a(n,k)\eq &4k^2-8kn-14k^2n-14k-4n+13k^2-11kn-26+39k+4n+26
\\=&k^2(17-14n)+k(25-19n)\pmod{n^2-1},
\endalign$$
and $(\pm n-1)\mid k\bi{\pm n-1}k$ since $k\bi{\pm n-1}k=(\pm n-1)\bi{\pm n-2}{k-1}$ if $k>0$.
So
$$\bi{n-1}k\bi{-n-1}ka(n,k)\eq \bi{n-1}k\bi{-n-1}kk(25-19n)\pmod{n^2-1}.$$
If $2\nmid n$, then $n\pm1$ and $25-19n$ are all even, hence both $2(n-1)$ and $2(n+1)$
divide $\bi{n-1}k\bi{-n-1}ka(n,k)$. If $n$ is even, then $(n-1,n+1)=(n-1,2)=1$ and hence
$n^2-1$ coincides with the least common multiple $[n-1,n+1]$ of $n-1$ and $n+1$.
Note that when $n$ is odd we have $(2,n-1)=2$ and
$$[2(n-1),2(n+1)]=\f{2(n-1)2(n+1)}{(2(n-1),2(n+1))}=\f{4(n^2-1)}{2(n-1,2)}=n^2-1.$$
Therefore $n^2-1\mid \bi{n-1}k\bi{-n-1}ka(n,k)$ no matter $n$ is odd or even.
This concludes the proof. \qed

\proclaim{Lemma 3.4} Let $a,b\in\N$ with $a+b$ even, and let $n\in\Z^+$. Then
$$2n\,\big|\,\sum_{k=0}^{n-1}\bi{n-1}k^a\bi{-n-1}k^b\bi{2k}k(k+2)3^{n-1-k}.\tag3.6$$
\endproclaim
\Proof. Let $f(k)=k\bi{2k-1}k3^{n-k}$ for $k=0,\ldots,n$. For each $k=0,\ldots,n-1$, we clearly have
$$\align\Delta f(k)=&f(k+1)-f(k)=(k+1)\bi{2k+1}{k+1}3^{n-k-1}-k\bi{2k-1}k3^{n-k}
\\=&(2k+1)\bi{2k}k3^{n-k-1}-3k\bi{2k-1}k3^{n-1-k}=\f{k+2}2\bi{2k}k3^{n-1-k}.
\endalign$$
Thus, by [S18a, Theorem 4.1] we get
$$\align&\sum_{k=0}^{n-1}\bi{n-1}k^a\bi{-n-1}k^b\f{k+2}2\bi{2k}k3^{n-1-k}
\\=&\sum_{k=0}^{n-1}\bi{n-1}k^a\bi{-n-1}k^b\Delta f(k)\eq0\pmod n
\endalign$$
and hence (3.6) holds. \qed

\medskip
\noindent {\it Proof of Theorem 1.2}. Since $(n,n^2-1)=1$, by Lemmas 3.2 and 3.3 it suffices to show that
$$\sum_{k=0}^{n-1}\bi{n-1}{k}\bi{-n-1}kC_k3^{n-1-k}a(n,k)\eq0\pmod n.$$

For each $k=0,\ldots,n-1$, clearly
$$a(n,k)\eq 13k^2+39k+26=13(k+1)(k+2)\pmod n.$$
So
$$\align &\sum_{k=0}^{n-1}\bi{n-1}{k}\bi{-n-1}kC_k3^{n-1-k}a(n,k)
\\\eq&13\sum_{k=0}^{n-1}\bi{n-1}{k}\bi{-n-1}k\bi{2k}k(k+2)3^{n-1-k}\eq0\pmod n.
\endalign$$
with the help of Lemma 3.4. This completes the proof. \qed

\heading{4. Proof of Theorem 1.3}\endheading

\proclaim{Lemma 4.1} Let $b,c\in\Z$ and $d=b^2-4c$. For any $n\in\Z^+$ we have
$$nT_n(b,c)T_{n-1}(b,c)=b\sum_{j=0}^{n-1}(n-j)\bi{n+j}{2j}\bi{2j}j^2c^jd^{n-1-j}.\tag4.1$$
\endproclaim
\Proof. In view of Lemma 3.1,
$$\align nT_n(b,c)T_{n-1}(b,c)=&b\sum_{k=0}^{n-1}(2k+1)\sum_{j=0}^k\bi{k+j}{2j}\bi{2j}j^2c^jd^{k-j}(-d)^{n-1-k}
\\=&b\sum_{j=0}^{n-1}\bi{2j}j^2c^jd^{n-1-j}\sum_{k=j}^{n-1}(-1)^{n-1-k}(2k+1)\bi{k+j}{2j}.
\endalign$$
For each $j\in\N$, by induction we have
$$\sum_{k=j}^{m-1}(-1)^{m-1-k}(2k+1)\bi{k+j}{2j}=(m-j)\bi{m+j}{2j}\ \ \t{for all}\ m=j+1,j+2,\ldots.\tag4.2$$
Thus
$$nT_n(b,c)T_{n-1}(b,c)=b\sum_{j=0}^{n-1}\bi{2j}j^2c^jd^{n-1-j}(n-j)\bi{n+j}{2j}$$
and hence (4.1) holds. \qed

\proclaim{Lemma 4.2} For any $k,n\in\Z^+$ with $k\ls n$, we have
$$\f{n(n+1)(n+2)}{(2,n)}\,\bigg|\, (n+k+1)\bi{n+k}k\bi{n+1}{k+1}\bi{2k}{k+1}.\tag4.3$$
\endproclaim
\Proof. Clearly,
$$\align &(n+k+1)\bi{n+k}k\bi{n+1}{k+1}\bi{2k}{k+1}
\\=&(n+k+1)\bi{n+k}k\f{n+1}{k+1}\bi nkkC_k
\\=&(n+1)\bi{n+k+1}{k+1}n\bi{n-1}{n-k}C_k,
\endalign$$
and also
$$(n+k+1)\bi{n+k}k\bi{2k}{k+1}
\eq(k-1)(-1)^k\bi{-n-1}kkC_k\eq0\pmod{n+2}$$
since
$$k(k-1)\bi{-n-1}k=(-n-1)(-n-2)\bi{-n-3}{k-2}\quad  \t{if}\ k>1.$$
Thus
$$[n(n+1),n+2]\,\big|\,(n+k+1)\bi{n+k}k\bi{n+1}{k+1}\bi{2k}{k+1}.$$
Note that
$$[n(n+1),n+2]=\f{n(n+1)(n+2)}{(n(n+1),n+2)}=\f{n(n+1)(n+2)}{(2,n)}.$$
So we have (4.3). \qed

\proclaim{Lemma 4.3} For any $n\in\N$ we have
$$6\bi{2n}n\eq0\pmod{n+2}.\tag4.4$$
\endproclaim
\Proof. Observe that
$$\bi{2n+2}{n+1}=2\bi{2n+1}n=\f{2(2n+1)}{n+1}\bi{2n}n$$
and hence
$$2(2n+1)\bi{2n}n=(n+1)\bi{2n+2}{n+1}=(n+1)(n+2)C_{n+1}.$$
Thus
$$\f{n+2}{(n+2,2n+1)}\,\big|\,\f{2n+1}{(n+2,2n+1)}2\bi{2n}n$$
and hence
$$\f{n+2}{(n+2,2n+1)}\,\big|\,2\bi{2n}n.\tag4.5$$
Since $(n+2,2n+1)=(n+2,2(n+2)-3)=(n+2,3)$ divides $3$, we obtain (4.4) from (4.5). \qed

As in [S18b], for $n\in\Z^+$ we define
$$w_n(x):=\sum_{k=1}^nw(n,k)x^{k-1}\ \t{with}\ w(n,k)=\f1k\bi{n-1}{k-1}\bi{n+k}{k-1}\in\Z.$$

\proclaim{Lemma 4.4} For any integers $n\gs k\gs1$, we have
$$w(n,k)=\sum_{j=1}^k\bi{n-j}{k-j}N(n,j)\tag4.6$$
and
$$N(n,k)=\sum_{j=1}^k\bi{n-j}{k-j}(-1)^{k-j}w(n,j).\tag4.7$$
\endproclaim
\Proof. We first prove (4.7). Observe that
$$\align\sum_{j=1}^k\bi{n-j}{k-1}(-1)^{k-j}w(n,j)
=&\sum_{j=1}^k\bi{n-j}{k-j}\f{(-1)^{k-j}}n\bi nj\bi{n+j}{j-1}
\\=&\f{(-1)^{k-1}}n\bi nk\sum_{j=1}^k\bi k{k-j}\bi{-n-2}{j-1}.
\endalign$$
Thus, with the help of the Chu-Vandermonde identity (cf. [G, (3.1)]), we get
$$\sum_{j=1}^k\bi{n-j}{k-1}(-1)^{k-j}w(n,j)=\f{(-1)^{k-1}}n\bi nk\bi{k-n-2}{k-1}=N(n,k).$$
This proves (4.7).

In view of (4.7), we have
$$\align \sum_{j=1}^k\bi{n-j}{k-j}N(n,j)=&\sum_{j=1}^k\bi{n-j}{k-j}\sum_{i=1}^j\bi{n-i}{j-i}(-1)^{j-i}w(n,i)
\\=&\sum_{i=1}^kw(n,i)\bi{n-i}{k-i}\sum_{j=i}^k\bi{k-i}{j-i}(-1)^{j-i}=w(n,k).
\endalign$$
So (4.6) also holds. This ends the proof. \qed

\proclaim{Lemma 4.5} For any $n\in\Z^+$ we have
$$w_n(x)=s_n(x).\tag4.8$$
\endproclaim
\Proof. With the aid of (4.7), we get
$$\align s_n(x)=&\sum_{k=1}^nN(n,k)x^{k-1}(x+1)^{n-k}
\\=&\sum_{k=1}^n\sum_{j=1}^k\bi{n-j}{k-j}(-1)^{k-j}w(n,j)x^{k-1}(x+1)^{n-k}
\\=&\sum_{j=1}^nw(n,j)x^{n-1}\sum_{k=j}^n\bi{n-j}{k-j}(-1)^{k-j}\l(1+\f1x\r)^{n-j-(k-j)}
\\=&\sum_{j=1}^nw(n,j)x^{n-1}\l(1+\f1x-1\r)^{n-j}=w_n(x).
\endalign$$
This concludes the proof. \qed

\proclaim{Lemma 4.6} For any $n\in\Z^+$ we have the new identity
$$(2x+1)\sum_{k=1}^n k(k+1)(2k+1)(-1)^{n-k}w_k(x)^2=n(n+1)(n+2)w_n(x)w_{n+1}(x).\tag4.9$$
\endproclaim
\Proof. In the case $n=1$, both sides of (4.9) are equal to $6(2x+1)$.

Now assume that $(4.9)$ holds for a fixed positive integer $n$.
Applying the Zeilberger algorithm (cf. [PWZ, pp.\,101-119]) via {\tt Mathematica 9} we find that
$$(n+3)w_{n+2}(x)=(2x+1)(2n+3)w_{n+1}(x)-nw_n(x).$$
Thus
$$\align &(2x+1)\sum_{k=1}^{n+1}k(k+1)(2k+1)(-1)^{n+1-k}w_k(x)^2
\\=&(2x+1)(n+1)(n+2)(2n+3)w_{n+1}(x)^2-(2x+1)\sum_{k=1}^{n}k(k+1)(2k+1)(-1)^{n-k}w_k(x)^2
\\=&(2x+1)(n+1)(n+2)(2n+3)w_{n+1}(x)^2-n(n+1)(n+2)w_n(x)w_{n+1}(x)
\\=&(n+1)(n+2)w_{n+1}(x)((2x+1)(2n+3)w_{n+1}(x)-nw_n(x))
\\=&(n+1)(n+2)(n+3)w_{n+1}(x)w_{n+2}(x).
\endalign$$

In view of the above, by induction, (4.9) holds for each $n\in\Z^+$. \qed

\medskip
\noindent{\it Proof of Theorem} 1.3. (i) Let $\da\in\{0,1\}$. In light of Lemma 4.1,
$$\align &\sum_{k=1}^nk^{2\da+1}T_k(b,c)T_{k-1}(b,c)d^{n-k}
\\=&\sum_{k=1}^nk^{2\da}b\sum_{j=0}^{k-1}(k-j)\bi{k+j}{2j}\bi{2j}j^2c^jd^{k-1-j}d^{n-k}
\\=&b\sum_{j=0}^{n-1}\bi{2j}j^2c^jd^{n-1-j}\sum_{k=j+1}^nk^{2\da}(k-j)\bi{k+j}{2j}.
\endalign$$
By induction, for each $j\in\N$, we have
$$\sum_{k=j+1}^mk^{2\da}(k-j)\bi{k+j}{2j}=\f{m^\da(m+1)^\da}2\cdot\f{(m-j)(m+j+1)}{j+\da+1}\bi{m+j}{2j}\tag4.10$$
for every $m =j+1,j+2,\ldots$. Therefore,
$$\align &\sum_{k=1}^nk^{2\da+1}T_k(b,c)T_{k-1}(b,c)d^{n-k}
\\=&b\f{n^\da(n+1)^\da}2\sum_{j=0}^{n-1}\bi{2j}j^2c^jd^{n-1-j}\f{(n-j)(n+j+1)}{j+\da+1}\bi{n+j}{2j}
\\=&\f b2(n(n+1))^\da\sum_{j=0}^{n-1}\f{\bi{2j}j}{j+\da+1}c^jd^{n-1-j}(n-j)(n+j+1)\bi nj\bi{n+j}j
\endalign$$
and hence
$$\aligned
&\sum_{k=1}^nk^{2\da+1}T_k(b,c)T_{k-1}(b,c)d^{n-k}
\\=&\f b2(n(n+1))^{\da+1}\sum_{j=0}^{n-1}\bi{n-1}j\bi{n+j+1}j\f{\bi{2j}j}{j+\da+1}c^jd^{n-1-j}.
\endaligned\tag4.11$$

In the case $\da=0$, (4.11) yields (1.4) since $\bi{2j}j/(j+1)=C_j\in\Z$.
By Lemma 4.3 and (4.11) with $\da=1$, we immediately obtain (1.5).

(ii) By induction, for each $j\in\N$ we have
$$\sum_{k=j}^m(2k+1)\bi{k+j}{2j}=\f{(m+1)(m+j+1)}{j+1}\bi{m+j}{2j}\quad\t{for all}\ m=j,j+1,\ldots.\tag4.12$$
In view of this and (2.4), we have
$$\align&\sum_{k=1}^nk(k+1)(2k+1)s_k(x)^2
\\=&\sum_{k=1}^n(2k+1)\sum_{j=1}^k\bi{k+j}{2j}\bi{2j}j\bi{2j}{j+1}(x(x+1))^{j-1}
\\=&\sum_{j=1}^n\bi{2j}j\bi{2j}{j+1}(x(x+1))^{j-1}\sum_{k=j}^n(2k+1)\bi{k+j}{2j}
\\=&\sum_{j=1}^n\bi{2j}j\bi{2j}{j+1}(x(x+1))^{j-1}\f{(n+1)(n+j+1)}{j+1}\bi{n+j}{2j}
\\=&\sum_{j=1}^n\bi{2j}{j+1}(x(x+1))^{j-1}\f{(n+1)(n+j+1)}{j+1}\bi nj\bi{n+j}j
\endalign$$
and hence
$$\sum_{k=1}^nk(k+1)(2k+1)s_k(x)^2=\sum_{k=1}^n(n+k+1)\bi{n+1}{k+1}\bi{n+k}k\bi{2k}{k+1}(x(x+1))^{k-1}.\tag4.13$$

Let $x=(b/\sqrt d-1)/2$. Then $x(x+1)=c/d$. In view of Lemma 2.1(ii) and (4.13), we have
$$\align&\sum_{k=0}^{n-1}(k+1)(k+2)(2k+3)M_k(b,c)^2d^{n-1-k}
\\=&\sum_{k=0}^{n-1}(k+1)(k+2)(2k+3)d^ks_{k+1}(x)^2d^{n-1-k}
\\=&d^{n-1}\sum_{k=1}^nk(k+1)(2k+1)s_k(x)^2
\\=&\sum_{k=1}^n(n+k+1)\bi{n+1}{k+1}\bi{n+k}k\bi{2k}{k+1}c^{k-1}d^{n-k}.
\endalign$$
Combining this with Lemma 4.2, we get the desired (1.6).

In light of Lemma 2.1(ii) and Lemmas 4.5-4.6, we have
$$\align &\sum_{k=0}^{n-1}(k+1)(k+2)(2k+3)M_k(b,c)^2(-d)^{n-1-k}
\\=&\sum_{k=0}^{n-1}(k+1)(k+2)(2k+3)d^ks_{k+1}(x)^2(-d)^{n-1-k}
\\=&d^{n-1}\sum_{k=1}^nk(k+1)(2k+1)(-1)^{n-k}w_k(x)^2
\\=&n(n+1)(n+2)d^{n-1}\f{s_n(x)s_{n+1}(x)}{2x+1}
\\=&n(n+1)(n+2)d^{n-1}\f{M_{n-1}(b,c)}{\sqrt d^{n-1}}\cdot\f{M_n(b,c)}{\sqrt d^n}\cdot\f{\sqrt d}b
\\=&n(n+1)(n+2)\f{M_n(b,c)M_{n-1}(b,c)}b.
\endalign$$
If $2\nmid n$ then $b\mid M_n(b,c)$; if $2\mid n$ then $2\nmid n-1$ and $b\mid M_{n-1}(b,c)$.
So $b$ divides $M_n(b,c)M_{n-1}(b,c)$. Therefore (1.7) holds.

The proof of Theorem 1.3 is now complete. \qed

\heading{5. Some open problems}\endheading

Clearly,
$$\f{\bi{2k}k}{2k-1}=\f{2}{2k-1}\bi{2k-1}k= \f 2k\bi{2k-2}{k-1}=2C_{k-1}\ \t{for}\ k\in\Z^+,$$
and thus $2k-1\mid\bi{2k}{k}$ for all $k\in\N$. Motivated by this we introduce a new kind of numbers
$$W_n:=\sum_{k=0}^{\lfloor n/2\rfloor}\bi n{2k}\f{\bi{2k}k}{2k-1}\ \ (n=0,1,2,\ldots)\tag5.1$$
which are analogues of the Motzkin numbers. The values of $W_0,W_1,\ldots,W_{12}$ are as follows:
$$-1,\ -1,\ 1,\ 5,\ 13,\ 29,\ 63,\ 139,\ 317,\ 749,\ 1827,\ 4575,\ 11699.$$
 Applying the Zeilberger algorithm (cf. [PWZ, pp.\,101-119]) via {\tt Mathematica 9}, we obtain the recurrence
 $$(n+3)W_{n+3} = (3n+7)W_{n+2}+ (n-5)W_{n+1} - 3(n+1)W_n\ (n=0,1,2,\ldots).$$
 For this new kind of numbers, we have the following conjecture similar to Theorem 1.1.

 \proclaim{Conjecture 5.1} {\rm (i)} For any $n\in\Z^+$ we have
 $$\sum_{k=0}^{n-1}(8k+9)W_k^2\eq n\pmod{2n}.\tag5.2$$
 Also, for any odd prime $p$ we have
 $$\f1p\sum_{k=0}^{p-1}(8k+9)W_k^2\eq 24+10\l(\f{-1}p\r)-9\l(\f p3\r)-18\l(\f 3p\r)\pmod{p}.\tag5.3$$
 
 {\rm (ii)} For any prime $p>3$ and positive integer $n$, the number
 $$\f{\sum_{k=0}^{pn-1}W_k^2-2(\sum_{k=0}^{n-1}T_k)^2}{pn}$$
 is always a $p$-adic integer.
 \endproclaim
 \Remark\ 5.1. We also guess that the sequence $(W_{n+1}/W_n)_{n\gs5}$ is strictly increasing to the limit $3$ and the sequence $(\!\root{n+1}\of{W_{n+1}}/\root n\of{W_n})_{n\gs9}$
 is strictly decreasing to the limit $1$.
\medskip

For $h,n\in\Z^+$, we define
$$w_n^{(h)}(x):=\sum_{k=1}^nw(n,k)^hx^{k-1}.$$

 \proclaim{Conjecture 5.2} Let $h,m,n\in\Z^+$. Then
 $$\f{(2,n)}{n(n+1)(n+2)}\sum_{k=1}^nk(k+1)(2k+1)w_k^{(h)}(x)^m\in\Z[x].\tag5.4$$
 Also,
 $$\f{(2,m-1,n)}{n(n+1)(n+2)}\sum_{k=1}^n(-1)^kk(k+1)(2k+1)w_k(x)^m\in\Z[x],\tag5.5$$
 and $$\f1{n(n+1)(n+2)}\sum_{k=1}^n(-1)^kk(k+1)(2k+1)w_k^{(h)}(x)^m\in\Z[x]\ \ \t{for}\ h>1.\tag5.6$$
\endproclaim

\Remark\ 5.2. Fix $n\in\Z^+$. By combining (4.13) with Lemma 4.2, we obtain
$$\f{(2,n)}{n(n+1)(n+2)}\sum_{k=1}^nk(k+1)(2k+1)s_k(x)^2\in\Z[x(x+1)].\tag5.7$$
As $s_k(x)=w_k(x)$ for all $k\in\Z^+$ (by Lemma 4.5), this implies (5.4) with $h=1$ and $m=2$.
Since $w_{2j}(x)/(2x+1)\in\Z[x]$ for all $j\in\Z^+$ (cf. [S18b, Section 4]), (5.5) with $m=2$ follows from (4.9).
\medskip

For $h\in\Z^+$ and $n\in\N$, we define
$$D_n^{(h)}(x):=\sum_{k=0}^n\bi nk^h\bi{n+k}k^hx^k
\ \ \t{and}\ \ S_n^{(h)}(x):=\sum_{k=0}^n\bi {n+k}{2k}^hC_k^hx^{k}.$$
Note that $S_n^{(1)}(x)=S_n(x)$ for all $n\in\N$.

\proclaim{Conjecture 5.3} Let $h,m,n\in\Z^+$. 

{\rm (i)} We have
 $$\f{(2,n)}{n(n+1)(n+2)}\sum_{k=1}^nk(k+1)(2k+1)S_k^{(h)}(x)^m\in\Z[x]\tag5.8$$
 and
 $$\f{(2,m-1,n)}{n(n+1)(n+2)}\sum_{k=1}^n(-1)^kk(k+1)(2k+1)S_k^{(h)}(x)^m\in\Z[x].\tag5.9$$
 
{\rm (ii)} We have
 $$\f{(2,n)}{n(n+1)(n+2)}\sum_{k=1}^nk(k+1)(2k+1)D_k^{(h)}(x)^m\in\Z[x]$$
 and
 $$\f{(2,hm-1,n)}{n(n+1)(n+2)}\sum_{k=1}^n(-1)^kk(k+1)(2k+1)D_k^{(h)}(x)^m\in\Z[x].$$
\endproclaim

\Remark\ 5.3. Fix $n\in\Z^+$. As $S_k(x)=(x+1)s_k(x)=(x+1)w_k(x)$ for all $k\in\Z^+$,
(5.8) and (5.9) with $h=1$ and $m=2$ do hold in view of Remark 5.2.
We also conjecture that
$$\f2{3n(n+1)}\sum_{k=1}^n(-1)^{n-k}k^2D_kD_{k-1}\ \ \t{and}\ \ 
\f1n\sum_{k=1}^n(-1)^{n-k}(4k^2+2k-1)D_{k-1}s_k$$
are positive odd integers.
\medskip

\proclaim{Conjecture 5.4} {\rm (i)} For any $h,m,n\in\Z^+$ we have
$$\f{2(2,n)}{n(n+1)(n+2)}\sum_{k=1}^nk(k+1)(k+2)(w_k^{(h)}(x)w_{k+1}^{(h)}(x))^m\in\Z[x]\tag5.10$$

{\rm (ii)}
For any $m,n\in\Z^+$ we have
$$\f{2(2,n)}{n(n+1)(n+2)(2x+1)^m}\sum_{k=1}^nk(k+1)(k+2)(w_k(x)w_{k+1}(x))^m\in\Z[x].\tag5.11$$
If $n\in\Z^+$ is even, then
$$\f{4}{n(n+1)(n+2)(2x+1)^3}\sum_{k=1}^nk(k+1)(k+2)w_k(x)w_{k+1}(x)\in\Z[x].\tag5.12$$
\endproclaim
\Remark\ 5.4. Recall that $w_{2j}(x)/(2x+1)\in\Z[x]$ for all $j\in\Z^+$ (by [S18b, Section 4]).
\medskip

\Ack. The author would like to thank Prof. Qing-Hu Hou
and the anonymous referee for their helpful comments.

\enddocument